# Adaptation Implies Internal Model

Eduardo D. Sontag
Department of Mathematics
Rutgers University
New Brunswick, NJ 08903
sontag@control.rutgers.edu
http://www.math.rutgers.edu/~sontag/

**Abstract**

This note provides a simple result showing, under suitable technical assumptions, that if a system $\Sigma$ adapts to a class of external signals $\mathcal{U}$, then $\Sigma$ must necessarily contain a subsystem which is capable of generating all the signals in $\mathcal{U}$. It is not assumed that regulation is robust, nor is there a prior requirement for the system to be partitioned into separate plant and controller components.

## 1   Introduction

Suppose that one knows that a certain system $\Sigma$ adapts to (or "regulates against") all those external input signals $u$ which belong to a predetermined class $\mathcal{U}$ of time-functions. (Input signals $u$ are often thought of as disturbances to be rejected or signals to be tracked, depending on the application.) In this context, adaptation means that a certain quantity $y(t)$ associated to the system, called its output (also called a regulated variable or an error) has the property that $y(t) \to 0$ as $t \to \infty$ whenever the system is subject to an input signal from the class $\mathcal{U}$ (Figure 1). The *internal model principle (IMP)* states, roughly, that the system $\Sigma$ necessarily must contain a subsystem $\Sigma_{\text{IM}}$ which can itself generate all disturbances in the class $\mathcal{U}$. The terminology arises when thinking of $\Sigma_{\text{IM}}$ as a "model" of a system which generates the external signals.

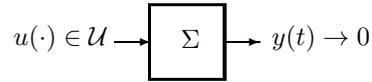

Figure 1: Given System, Regulated Output $y(t)$ when Inputs in $\mathcal{U}$

For example, if $y(t) \to 0$ as $t \to \infty$ whenever the system is subject to any external constant signal (i.e., the class $\mathcal{U}$ consists of all constant functions), then the system $\Sigma$ must contain a subsystem $\Sigma_{\text{IM}}$ which generates all constant signals (typically an integrator, since constant signals are generated by the differential equation $\dot{u} = 0$). Of course, the choice of $y = 0$ as the "adaptation value" is merely a matter of convention; by means of a change of variables, one may always reduce a given regulation objective "$y(t) \to y_0$" where $y_0$ is some predetermined value, to the special case $y_0 = 0$.

In addition, the IMP specifies that, in an appropriate sense, the subsystem $\Sigma_{\text{IM}}$ must only have $y$ as its external input, receiving no other direct information from other parts of the system nor the input signal $u$. One intuitive interpretation is that $\Sigma_{\text{IM}}$ generates its "best guess" of the external input $u$ based on how far the output $y$ is from zero. Pictorially, if we have the situation shown in Figure 1, then there must be a decomposition of the system $\Sigma$ into two parts, as shown in Figure 2, where the system $\Sigma_{\text{IM}}$ (with $y \equiv 0$) is capable of reproducing all the functions in $\mathcal{U}$.

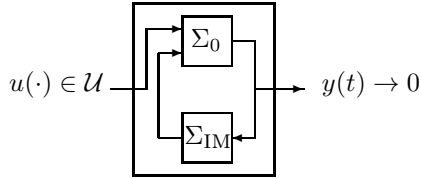

Figure 2: Decomposition of $\Sigma$ into $\Sigma_0$ and $\Sigma_{\text{IM}}$, the Latter Driven by $y(t)$

The internal model principle originates in the biological cybernetics literature. But, as with any "principle" in control theory (like dynamic programming, the maximum principle, etc.) and more generally in mathematics, the IMP is not a theorem but rather a "mold" for many possible theorems, each of which will hold under appropriate technical assumptions, and whose conclusions will depend upon the precise meaning of "class of external signals", "reproducing all functions", and so on.

The best known instance of an internal model *theorem* is due to Francis and Wonham, who in a series of beautiful and deep papers in the mid 1970s proved a theorem for linear systems which showed, in essence, that structurally stable or "robust" adaptation implies the existence of internal models. Partial generalizations of their work to nonlinear systems were later obtained by Wonham and Hepburn, see [1],[15],[16],[2]-[7]. The Francis/Wonham theory applies to systems $\Sigma$ which are already partitioned into a "plant" plus a "controller". The robustness assumption amounts to the requirement that the given controller should perform appropriately (in the sense that the regulation objective $y(t) \to 0$ is achieved) even when the plant subsystem – but most definitely not the controller subsystem – is arbitrarily perturbed. The conclusion is that the controller is driven by $y$ and incorporates a model of the external signals. That some additional condition – such as structural stability – must be imposed is obvious, since the system $\Sigma$ which simply outputs $y = 0$ for every possible input signal $u$ does not contain any subsystem generating the signals $u$. We will impose instead a condition which amounts to a "signal detection" property: the output must reflect sudden changes in the input.

Recent work in molecular biology, cf. [17], has suggested that the IMP could help guide experimentalists and modelers: if certain characteristics of a system adapt to signals in a given class (in all the examples so far, constant inputs, such as for instance $y(t) =$ the relative "activity" of enzymes controlling motors in *E.coli* chemotaxis, with respect to $u(t) =$ concentration of extracellular ligand) then the IMP could, in principle, help distinguish among mathematical models which do or do not contain internal models.

With a view toward such biological applications, it is desirable to have available a theorem which (a) applies to nonlinear systems $\Sigma$, at least under reasonable technical assumptions, and (b) does not require the system $\Sigma$ to be split between "plant" and "controller" subsystems, nor requires structural stability (robustness) in the sense of the Francis/Wonham theory (which would imply, in the case of the *E.coli* motor control network, that the system should perfectly adapt even if there are arbitrary direct connections between the external ligands and the motor signals, a matter which seems difficult to check experimentally). We present one very elementary and self-contained such result in this note. It basically just picks from and "repackages" some of the basic concepts and techniques developed by previous researchers for the same problems, in particular: the use of differential geometric techniques and "output-zeroing" sets ([15], [2]-[7], [10]), dynamical systems notions like omega-limit sets ([2],[6],[10]), and system decompositions motivated by the Center Manifold Theorem ([2],[5],[6],[10]). Isidori's excellent textbooks [11],[12] should be consulted for a far deeper discussion of many of the issues raised here.

Precise mathematical definitions are provided in Section 2. On the other hand, since the linear version of the result is very easy to explain, we sketch that case first. (The discussion assumes some familiarity with frequency domain techniques, and may be skipped without loss of continuity.)

## 1.1 Linear Case

Let us denote by $S$ the *transfer function* of the system $\Sigma$: if $y$ is the output produced when $\Sigma$ starts at the zero initial state and is fed input $u$, then the relation $\hat{y}(s) = S(s)\hat{u}(s)$ holds between the Laplace transforms $\hat{y}(s)$, $\hat{u}(s)$ of the output and input. One expresses $S(s) = \frac{p(s)}{q(s)}$ as the quotient of two relatively prime polynomials, with the degree of $p$ less than the degree of $q$. (An equivalent discussion using differential operators instead of Laplace transforms is also possible, see e.g. Section 6.7 in [14].) The first observation, a well-known fact in systems theory, is that the zeroes of $p$ can be viewed, alternatively, as poles of a feedback subsystem. To see this, we assume that $p$ is not identically zero, and divide the polynomial $q$ by $p$, obtaining $q = ap + b$, where $b$ is some polynomial of degree less than $p$. Now, as the algebraic equality $y = \frac{p}{q}u$ is equivalent to $y = \frac{1}{a}(u - \frac{b}{p}y)$, we conclude that the system $\Sigma$ can be decomposed as in Figure 3. For example, if $s = 0$ is a zero of $S$ (that is, 0 is a root of $p$), which amounts to the property that constant signals get differentiated by $\Sigma$ (the "DC gain" of $\Sigma$ is zero), then the factor $1/s$ appears in the feedback box $b/p$, and can be interpreted as an integrator of the output $y$.

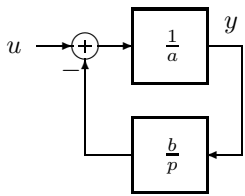

Figure 3: System Equivalent to $\Sigma$: Closed-loop Zeros are Feedback Poles

We will show that, in general, the subsystem with transfer function $\frac{b}{p}$ models all inputs which $\Sigma$ adapts to. Let us suppose that the class $\mathcal{U}$ of inputs can be described as the set of all possible solutions of a fixed linear differential equation

$$u^{(\ell)}(t) + b_1 u^{(\ell-1)}(t) + \ldots + b_{\ell-1} u'(t) + b_\ell u(t) = 0$$

for some integer $\ell$, and which has no stable modes. (Stable modes, giving components of $u$ which converge to zero, are less interesting, since they do not represent persistent disturbances.) We view these signals $u$ as the outputs of an "exosystem" $\Gamma$ which is obtained by rewriting the differential equation as a system of $\ell$ first order equations. Figure 4 shows a cascade consisting of the original system $\Sigma$ and the exosystem $\Gamma$ which generates the inputs in $\mathcal{U}$. (If, for example, $\mathcal{U}$ = constant inputs, then one would let $\Gamma$ be the system with equation $\dot{w} = 0$ and output $u = w$, and for each initial condition $w(0) = u^0$, one obtains a different constant output $u(t) \equiv u^0$.) The regulation objective is now simply that $y(t) \to 0$ for all possible initial conditions of the composite system, i.e. for all initial conditions of the original system $\Sigma$ and all initial conditions of the exosystem $\Gamma$, the latter corresponding to all possible inputs in $\mathcal{U}$ being fed to $\Sigma$.

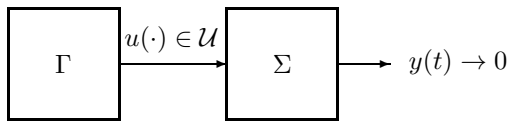

Figure 4: Exosystem and System in Cascade

Next, we reformulate this regulation property by adding an external input $v(\cdot)$ to the exosystem, and requiring now that $y(t) \to 0$ for all possible stable inputs ($v(t) \to 0$ as $t \to \infty$) but only when starting from the zero initial state. (Such replacements of initial states by stable forcing inputs – assuming natural controllability/observability conditions – are elementary exercises in linear systems theory, see e.g. the proof of Theorem 33 in [14].) In other words, we have now the situation illustrated by Figure 5. We denote by $G$ the transfer function of the exosystem $\Gamma$: $G = \frac{1}{\pi}$, where

$$\pi(s) = s^\ell + b_1 s^{\ell-1} + \ldots + b_{\ell-1} s + b_\ell\,.$$

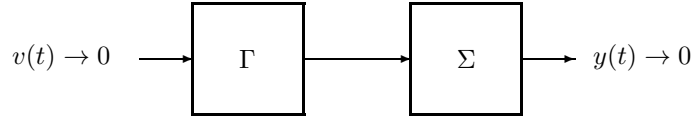

Figure 5: Exosystem and System, Forced by Stable Inputs

To see that the subsystem with transfer function $b/p$ includes an internal model of $\Gamma$, we argue as follows. The regulation property for the cascade in Figure 5 means that the product rational function $GS$ is stable (all poles have negative real parts), while the assumption that $G$ had no stable modes means that all the poles of $G$ (i.e, the roots of the polynomial $\pi$) have nonnegative real parts. Therefore, these poles must be canceled in the product $GS$; in other words, $S$ *must have among its zeroes all the poles of* $G$, so that we can write $p = \pi p_0$ for some polynomial $p_0$. Thus $b/p = b/(\pi p_0)$. One may now factor $b = b_1 b_2$ in such a way that the degree of $b_2$ is less than the degree of $\pi$, so that $b/p = (b_1/p_0)(b_2/\pi)$ and now the system with transfer function $b/p$ can be written itself in the cascade form in Figure 6. The subsystem with transfer function $b_2/\pi$ generates all the inputs in $\mathcal{U}$, since one may write a set of differential equations for it which is exactly the same as for the exosystem $\Gamma$, changing only the output mapping ("controller form" realization).

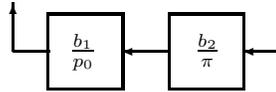

Figure 6: Decomposition of $\frac{b}{p}$

Since the tools of transfer functions are not available for nonlinear systems, a different approach is required in general.

## 2  Definitions and Statement of Result

We consider single-input single-output systems $S$, affine in inputs:

$$\dot{x}(t) = f(x(t)) + u(t)g(x(t))\,,\ y(t) = h(x(t)) \tag{1}$$

(dot indicates derivative with respect to time, and the arguments $t$ will be omitted from now on; see [14] for general definitions and properties of systems with inputs). Here $x(t)$, $u(t)$, and $y(t)$ represent the state, input, and output at time $t$, $f$ and $g$ are smooth vector fields on $\mathbb{R}^n$ ($n$ is the dimension of the system), $h$ is a scalar smooth function $\mathbb{R}^n \to \mathbb{R}$, and $f(0) = h(0) = 0$. (Several assumptions on $f$ and $g$ will be made later.) A special case is that of *linear* systems

$$\dot{x} = Ax + ub\,,\ y = cx \tag{2}$$

where $A$ is an $n \times n$ matrix, $b$ is a column $n$-vector, and $c$ is a row $n$-vector.

Suppose given a class $\mathcal{U}$ of functions $[0, \infty) \to \mathbb{R}$ (such as for example the set of all constant functions). We say that $\Sigma$ *adapts* to inputs in $\mathcal{U}$ (a more appropriate technical control-theoretic term would be "asymptotically rejects disturbances in $\mathcal{U}$") if the following property holds: for each $u \in \mathcal{U}$ and each initial state $x^0 \in \mathbb{R}^n$, the solution of (1) with initial condition $x(0) = x^0$ exists for all $t \geq 0$ and is bounded, and the corresponding output $y(t) = h(x(t))$ converges to zero as $t \to \infty$.

We will say that *S contains an output-driven internal model of $\mathcal{U}$* if there is a change of coordinates which brings the equations (1) into the following block form:

$$\begin{aligned} \dot{z}_1 &= f_1(z_1, z_2) + u g_1(z_1, z_2) \\ \dot{z}_2 &= f_2(y, z_2) \\ y &= \kappa(z_1) \end{aligned} \quad (3)$$

(the subsystems with variables $z_1$ and $z_2$ correspond respectively to $\Sigma_0$ and $\Sigma_{\text{IM}}$ in Figure 2), and in addition the subsystem with state variables $z_2$ is capable of generating all functions in $\mathcal{U}$, meaning the following property: there is some scalar function $\varphi(z_2)$ so that, for each possible $u \in \mathcal{U}$, there is some solution of

$$\dot{z}_2 = f_2(0, z_2) \quad (4)$$

which satisfies $\varphi(z_2(t)) \equiv u(t)$.

The precise meaning of "change of coordinates" is as follows. There must exist an integer $r \leq n$, differentiable manifolds $Z_1$ and $Z_2$ of dimensions $r$ and $n-r$ respectively, a smooth function $\kappa : Z_1 \to \mathbb{R}$, vector fields $F$ and $G$ on $Z_1 \times Z_2$ which take the partitioned form

$$F = \begin{pmatrix} f_1(z_1, z_2) \\ f_2(\kappa(z_1), z_2) \end{pmatrix}, \quad G = \begin{pmatrix} g_1(z_1, z_2) \\ 0 \end{pmatrix}$$

and a diffeomorphism $\Phi : \mathbb{R}^n \to Z_1 \times Z_2$, such that

$$\Phi_*(x) f(x) = F(\Phi(x)), \quad \Phi_*(x) g(x) = G(\Phi(x)), \quad \kappa(\Phi_1(x)) = h(x)$$

for all $x \in \mathbb{R}^n$, where $\Phi_1$ is the $Z_1$-component of $\Phi$ and star indicates Jacobian.

Our result will hold under additional conditions on the vector fields defining the system. The first condition is the fundamental one from an intuitive point of view, namely that the system is able to detect changes in the input signal:

> **Assumption 1:** a uniform relative degree exists.

This means that there exists some positive integer $r$ such that

$$L_g L_f^k h \equiv 0 \quad \forall k < r - 1$$

and

$$L_g L_f^{r-1} h(x) \neq 0 \quad \forall x \in \mathbb{R}^n$$

where, as usual, $L_X h$ indicates the directional derivative ("Lie derivative") of a function $h$ along the direction of the vector field $X$, that is $(L_X h)(x) = \nabla h(x) \cdot X(x)$. The integer $r$, if it exists, is called the *relative degree* of $\Sigma$. It is possible to prove (see [11]) that when $r$ exists, necessarily $r \leq n$.

For a linear system (2), existence of a relative degree amounts to simply asking that $cA^i b$ is nonzero for some $i$, or equivalently that the transfer function $c(sI - A)^{-1} b$ is not identically zero. For general systems (1), the assumption is equivalent to the statement that the output derivatives $y^{(k)}(t)$ must be independent of the value of the input at time $t$, for all $k < r$, but that $y^{(r)}(t) = b(x(t)) + a(x(t))u(t)$ for some function $a(x)$ which is everywhere nonzero (so that the system can be "inverted" to obtain the instantaneous value $u(t)$ from instantaneous derivatives). See also [13] for a discussion of the characterization of $r$ in terms of smoothness of outputs when inputs are discontinuous (change detection).

The next two conditions are of a technical nature. They are automatically satisfied for linear systems. For nonlinear systems, we need such conditions in order to guarantee the existence of a change of variables exhibiting the system $\Sigma_{\rm IM}$. (Weaker conditions may be given, if one is merely interested in a local result, or if one is willing to accept a subsystem $\Sigma_{\rm IM}$ which is driven by not just $y$ but also several derivatives of $y$.) We are guided by conditions which appear in Isidori's [11].

Assuming that the degree is $r$, we introduce the following vector fields:

$$\widetilde{g}(x) \;=\; \frac{1}{L_g L_f^{r-1} h(x)} g(x), \quad \widetilde{f}(x) \;=\; f(x) - \left(L_f^r h(x)\right) \widetilde{g}(x), \quad \tau_i := {\rm ad}_{\widetilde{f}}^{i-1} \widetilde{g}, \; i = 1, \ldots r,$$

where ${\rm ad}_X$ is the operator ${\rm ad}_X Y = [X, Y] =$ Lie bracket of the vector fields $X$ and $Y$. Recall that a vector field $X$ is said to be *complete* if the solution of the initial value problem $\dot{x} = X(x)$, $x(0) = x^0$ is defined for all $t \in \mathbb{R}$, for any initial state $x^0$, and that two vector fields $X$ and $Y$ are said to *commute* if $[X, Y] = 0$. The assumptions are:

> **Assumption 2:** $\tau_i$ is complete, for $i = 1, \ldots, r$.

> **Assumption 3:** the vector fields $\tau_i$ commute with each other.

(For linear systems, the vector fields $\tau_i$ are all constant, so that they are indeed complete and pairwise commutative.)

Finally, we must define the allowed classes of inputs $\mathcal{U}$. As usual in control theory (see also the discussion in Section 1.1), we will assume that inputs are generated by *exosystems*. That is, there is given a system $\Gamma$:

$$\dot{w} = Q(w), \quad u = \theta(w) \tag{5}$$

(let us say evolving on some differentiable manifold, $Q$ a smooth vector field, and $\theta$ a real-valued smooth function, although far less than smoothness is needed) such that the input class $\mathcal{U}$ consists exactly of those inputs $u(t) = \theta(w(t))$, $t \geq 0$, for all possible solutions of $\dot{w} = Q(w)$. For example, if we are interested in constant signals, we pick $\dot{w} = 0$, $u = w$ and if we are interested in sinusoidals with frequency $\omega$ then we use $\dot{x}_1 = x_2$, $\dot{x}_2 = -\omega^2 x_1$, $u = x_1$. It is by now standard in nonlinear studies of necessary conditions for regulation to impose conditions on omega limits sets for trajectories of the exosystem, see [2],[6]; we will follow the approach in [10]-[11] and assume that the exosystem is Poisson-stable: for every state $w^0$, the solution $w(\cdot)$ of $\dot{w} = Q(w)$, $w(0) = w^0$ is defined for all $t > 0$ and it satisfies that $w^0$ is in the omega-limit set of $w$, that is, there is a sequence of times $t_i \to \infty$ such that the sequence $w(t_i)$ converges to $w^0$ as $t \to \infty$. This means that the exosystem is almost-periodic in the sense that trajectories keep returning to neighborhoods of the initial state.

This theorem is proved in Section 3:

**Theorem 1** *If Assumptions 1-3 hold and the system $\Sigma$ adapts to inputs in a class $\mathcal{U}$ generated by a Poisson-stable exosystem, then $S$ contains an output-driven internal model of $\mathcal{U}$.*

## 2.1 An Example

As an example, consider the model for *E.coli* chemotaxis adaptation to constant inputs given in [9], Section 2.2. Letting $x_1 = R$ and $x_2 = RL$ be the concentrations of unbound and bound receptors respectively, and taking the external ligand concentration $u = L$ as input, we have the following equations:

$$\begin{aligned} \dot{x}_1 &= a_1 - a_2 x_1 + a_3 x_2 - a_4 x_1 u \\ \dot{x}_2 &= a_5 - a_6 x_2 + a_4 x_1 u \end{aligned} \tag{6}$$

for suitable positive constants $a_1, \ldots, a_6$. In terms of vector fields,

$$f \;=\; \begin{pmatrix} a_1 - a_2 x_1 + a_3 x_2 \\ a_5 - a_6 x_2 \end{pmatrix}, \quad g \;=\; \begin{pmatrix} -a_4 x_1 \\ a_4 x_1 \end{pmatrix}$$

and, still as in [9], we take as output $y$ the difference between the total concentration of active receptors and a steady state level of this activity. In terms of the notations used here, and up to multiplication by a suitable constant, this amounts to the following choice:

$$h(x) = A_0 - A = [a_1 + a_5] - [a_2 x_1 + (a_6 - a_3)x_2].$$

We note that $L_g h = D x_1$, where $D = a_2 a_4 + (a_3 - a_6) a_4$. Except in the accidental case when this constant $D$ vanishes (in terms of the notations in [9], $D = k_{-1} I_T k_r (\alpha_1 - \alpha_2)$, so $D$ can only vanish if $\alpha_1 = \alpha_2$), we have that $L_g h(x) \neq 0$ for all $x$ ($x_1 > 0$, as it represents a concentration), and so it follows that $\Sigma$ has well-defined relative degree $r = 1$. Moreover, $\tau_1 = \widetilde{g}$ is a constant vector, so Assumptions 2 and 3 hold as well.

A minor technicality concerns the assumptions that our systems (1) evolve in all of Euclidean state space (not just $x_i > 0$) and that $f(0) = h(0) = 0$. However, this is just a matter of picking the right coordinates. Notice that $f$ vanishes at $x^0 = (x_1^0, x_2^0)$, where $x_1^0 = (a_1 a_6 + a_3 a_5)/(a_2 a_6)$ and $x_2^0 = a_5/a_6$, and $h$ vanishes at $x^0$ too. In order to fit into the general theory, one simply changes variables, mapping the positive orthant into all of $\mathbb{R}^2$ and $x^0$ into the origin by means of $x_i' = \ln x_i - \ln x_i^0$. (Of course, there is no need to actually perform the coordinate change, since conditions expressed in terms of Lie derivatives are covariant.)

Finally, letting $B := x_1 + x_2$ (as done in [9]), one obtains a system of equations in terms of the new variables $A$ and $B$, for which $\dot{B} = y$. This last equation represents an integrator (internal model of a system which produces constant inputs) driven by the output $y$. (Of course, there is no point in applying the theorem, since once that the model is given we can find the internal model explicitly.)

## 3 Proof of Theorem 1

Suppose that the system $\Sigma$ adapts to inputs in $\mathcal{U}$, which are produced by a Poisson stable exosystem $\Gamma$. We consider the interconnected system consisting of the cascade of $\Gamma$ and $\Sigma$, as shown in Figure 4, namely:

$$\begin{aligned} \dot{w} &= Q(w) \\ \dot{x} &= f(x) + \theta(w)g(x) \end{aligned} \quad (7)$$

and let $\mathcal{Z}$ denote the set consisting of those states $x$ of $\Sigma$ for which $h(x) = 0$ (the "output-zeroing" subset).

**Lemma 3.1** *For each $w^0$ there is some solution $\sigma = (w(\cdot), x(\cdot))$ of the composite system (7) such that $w(0) = w^0$ and $x(t) \in \mathcal{Z}$ for all $t \geq 0$.*

*Proof.* We start by picking an arbitrary solution $\sigma_0 = (w(\cdot), x(\cdot))$ of the composite system (7) such that $w(0) = w^0$, and let $\Omega = \Omega^+[\sigma_0]$ be the omega-limit set of this trajectory. We claim that, for each point $(\omega, \xi) \in \Omega$ (we partition coordinates into those for $\Gamma$ and $\Sigma$) it must be the case that $\xi \in \mathcal{Z}$. Indeed, by definition of $\Omega$ there is some sequence of times $t_i \to \infty$ such that $x(t_i) \to \xi$. Since $h(x(t_i)) \to 0$ because of the adaptation property and $h$ is continuous, it follows that $h(\xi) = 0$, as claimed. Next, we claim that there is some $x^0$ such that $(w^0, x^0) \in \Omega$. To see this, we first pick a sequence of times $t_i \to \infty$ such that $w(t_i) \to w^0$ (Poisson stability is used here); as $\{x(t_i)\}$ is bounded, we may pick a subsequence $t_{i_j}$ of the $t_i$ so that $x(t_{i_j}) \to x^0$ for some $x^0$, and this proves that $(w^0, x^0) \in \Omega$, as wanted.

Finally, we let $\sigma$ be the solution $\sigma = (w(\cdot), x(\cdot))$ of the composite system (7) for which $w(0) = w^0$ and $x(0) = x^0$, where $x^0$ is so that $(w^0, x^0) \in \Omega$. Omega-limit sets are invariant, so $\sigma(t) \in \Omega$ for all $t \geq 0$, and we already proved that this last property implies that $x(t) \in \mathcal{Z}$. ∎

Proposition 9.1.1 in [11] shows that there is a global diffeomorphism $\Phi$ so that, in the new coordinates, the system $\Sigma$ takes the form shown in Display (3). Moreover, the subsystem described by $z_1$ evolves in $\mathbb{R}^r$ and, using coordinates $z_1 = (\zeta_1, \ldots, \zeta_r)$, the equations for $z_1$ can be written as follows:

$$\begin{aligned} \dot{\zeta}_1 &= \zeta_2 \\ &\vdots \\ \dot{\zeta}_{r-1} &= \zeta_r \\ \dot{\zeta}_r &= b(z_1, z_2) + a(z_1, z_2)u \end{aligned} \qquad (8)$$

where the output is $y = \kappa(z_1) = \zeta_1$ and $a, b$ are smooth functions with $a(z) = L_g L_f^{r-1} h(\Phi^{-1}(z)) \neq 0$ for all $z$. We let

$$\varphi(z_2) := -\frac{b(0, z_2)}{a(0, z_2)}$$

and show that for each possible $u \in \mathcal{U}$ there is some solution of (4) which satisfies $\varphi(z_2(t)) \equiv u(t)$.

We pick $w^0$ such that $u(t) = \theta(w(t))$ and $w(0) = w^0$, and view the interconnection (7) of $\Gamma$ and $\Sigma$ in terms of the coordinate change given by $\Phi$ on $\Sigma$:

$$\begin{aligned} \dot{w} &= Q(w) \\ \dot{z}_1 &= f_1(z_1, z_2) + \theta(w) g_1(z_1, z_2) \\ \dot{z}_2 &= f_2(y, z_2). \end{aligned}$$

Lemma 3.1 gives us the existence of a solution $\sigma = (w(\cdot), z_1(\cdot), z_2(\cdot))$ such that $\theta(w(t)) = u(t)$ and $\zeta_1(t) \equiv 0$. Because of the form (8) of the $z_1$-subsystem, this implies that $z_1(t) \equiv 0$ and that $\dot{\zeta}_r(t) \equiv 0$. Thus, along the solution $\sigma$ one has $b(0, z_2(t)) + a(0, z_2(t))u(t) \equiv 0$, and this is precisely what we wished to prove. ∎

**A Remark on Subsystems**

We expressed our theorem in terms of the existence of solutions which reproduce all inputs. Under additional and stronger hypotheses, one could also obtain an actual embedding of the exosystem in the internal model $\Sigma_{\text{IM}}$. A full nonlinear version would involve abstract quotients of systems under suitable equivalence relations, and may follow along the lines of the work in [4] (based on [8]). However, the necessary steps are easy to understand and prove in the case of *linear* systems. We start by showing the following elementary fact from linear systems theory

**Lemma 3.2** Suppose given an observable linear system $\dot{w} = Qw$, $y = \theta w$ and another linear system $\dot{z}_2 = Fz_2 + Gy$, $u = \varphi z$, and assume that for each $w^0$ there is some $z^0$ such that $\varphi e^{tF} z^0 = \theta e^{tQ} w^0$ for all $t \geq 0$. Then, the matrix $F$ is similar to a matrix with this block structure:

$$\begin{pmatrix} Q & 0 & 0 \\ D & E & 0 \\ F & G & H \end{pmatrix}. \qquad (9)$$

*Proof.* We first assume that the pair $(F, \varphi)$ is observable, and claim that for each $w^0$ there is a *unique* $z^0$ such that $\varphi e^{tF} z^0 \equiv \theta e^{tQ} w^0$. This is because $\varphi e^{tF} z^0 \equiv \varphi e^{tF} z^1$ implies $z^0 = z^1$ (observability). So we can define a map $T: w^0 \mapsto z^0$. This map is one-to-one, by observability of the pair $(Q, \theta)$. It is also linear, since $\theta e^{tQ}(\alpha w^0 + w^1) = \alpha \theta e^{tQ} w^0 + \theta e^{tQ} w^1 = \alpha \varphi e^{tF} T w^0 + \varphi e^{tF} T w^1 = \varphi e^{tF}(\alpha T w^0 + T w^1)$ means that $\alpha w^0 + w^1 \mapsto \alpha T w^0 + T w^1$. It also satisfies $FT = TQ$, since taking derivatives in $\varphi e^{tF} T w^0 \equiv \theta e^{tQ} w^0$ gives $\varphi e^{tF} FT w^0 \equiv \theta e^{tQ} Q w^0$ which means that $Q w^0 \mapsto FT w^0$. Thus, on some invariant subspace (the range of $T$), $F$ can be written as $Q$, which means that we can write $F$ up to similarity in the form $\begin{pmatrix} Q & * \\ 0 & * \end{pmatrix}$. Since $F$ is similar to its transpose, and $Q$ is similar to its transpose, $F$ is also similar to a matrix in the form $\begin{pmatrix} Q & 0 \\ * & * \end{pmatrix}$. An observability decomposition ([14], Chapter 6) then reduces to the observable case. ∎

Without loss of generality, one may assume that linear exosystems are observable (there always exists an observable equivalent). We now apply Lemma 3.2 to the exosystem and the internal model $\Sigma_{\text{IM}}$, assumed linear. There results a change of variables for $\Sigma_{\text{IM}}$ so that, in the new variables, a subset $\zeta$ of the variables $z_2$ of $\Sigma_{\text{IM}}$, corresponding to the first block in (9), evolves according to an equation of the form $\dot\zeta = Q\zeta + by$, for a suitable vector $b$. This provides the desired embedding of the exosystem in the internal model.